\hyphenation{para-met-rize}

\def\Q{{\bf Q}}
\def\F{{\bf F}}
\def\Z{{\bf Z}}

\def\ref#1#2#3#4#5#6{ #1: #2. #3 {\bf #4}, #5 (#6)}

\def\card{\mathop{\rm card}}

\def\define{\smallskip{\noindent\bf Definition. }}
\def\which#1#2{\smallskip\global\advance \thmnum by 1 \edef#2{#1 \the\secnum .\the\thmnum}\noindent{\bf #2}\ \thinspace}
\def\example#1{\which{Example}{#1}}
\def\thm#1{\which{Theorem}{#1}}

\def\cor#1{\which{Corollary}{#1}}

\def\prop#1{\which{Proposition}{#1}}
\def\proof{\smallskip{\noindent\sl Proof.\ }}

\def\ref#1#2#3#4#5#6{ #1: #2. #3 {\bf #4}, #5 (#6)}

 4
\font\largebold=cmbx10 scaled \magstep 3
 2
\font\big=cmbx12

\def\Q{{\bf Q}}
\def\Z{{\bf Z}}

\def\isisom{\cong}
\def\chapter#1{\vfill\eject{\largebold\centerline{#1}\par}\vskip .2in}
\def\section#1{\vskip .05in \global\advance \secnum by 1 
{\big\centerline{{\ifnum\the\secnum>0 \the\secnum. \fi}#1}\par}\vskip .07in \global\thmnum=0}

\overfullrule=0pt
\def\mapsto{\to}

\def\Gal{\mathop{\rm Gal}}

\newcount\secnum
\secnum = 0
\newcount\thmnum
\thmnum = 0
\def\qed{ \ {\vbox{\hrule height 0pt depth .3pt\hbox{\vrule width .3pt height 6pt
\kern 3pt \vrule width .3pt}\hrule height 0pt depth .3pt}}\smallskip}

\def\subbar#1{{}_{#1} \bar}
\def\title#1{{\centered{\largebold #1}\medskip}}
\def\centered{\leftskip=0pt plus 4 em \rightskip=\leftskip
\parfillskip=0 pt \spaceskip=.3333 em \xspaceskip=.5em
\pretolerance=9999 \tolerance=9999 \parindent=0pt \hyphenpenalty=5000
\exhyphenpenalty=5000}
\def\today{\ifcase\month\or January\or February\or March\or April\or May\or
June\or July\or August\or September\or October\or November\or December\fi
\ \number\day, \number\year}
\newfam\frakfam
\font\fraktur=eufm10
\font\frakseven=eufm7
\font\frakfive=eufm5
\def\frak{\fam\frakfam\fraktur} \textfont\frakfam=\fraktur
\scriptfont\frakfam=\frakseven \scriptscriptfont\frakfam=\frakfive
\def\p{{\frak p}}
\def\m{{\frak m}}

\title{Deformations of Galois representations arising from degenerate
extensions}
\centerline{Adam Logan}
\smallskip
\centerline{\today}
\medskip
\section{Introduction}

This paper is inspired by that of Boston and Mazur [B-M], and work on this
problem was begun when the author was a graduate student of Barry Mazur
(supported by an NSF Graduate Fellowship).  In the paper [B-M], the authors
study the deformation theory of a certain type of $S_3$-extensions of $\Q$,
which they term {\sl neat}, and more specifically that of {\sl generic}
$S_3$-extensions, which satisfy an additional condition.  Restricting their
numerical study to one particular family of neat extensions, they note that
all such extensions seem to satisfy their genericity condition. 

Their principal result on generic $S_3$-extensions can be summarized as
follows:

\thm{\bmresult} ([B-M], prop.\ 13]) Let $L/\Q$ be a neat $S_3$-extension
for the prime $p$ (we will define this in Section 2, below).
The universal deformation ring of the natural representation of its Galois
group into $GL_2(\F_p)$ is isomorphic to $\Z_p[[T_1,T_2,T_3]]$.  If $L/\Q$
is generic, then:
\item{(a)}
The inertially reducible locus is composed of the union of two smooth
hypersurfaces in the universal deformation space.
\item{(b)}
The globally dihedral locus is equal to the inertially dihedral locus
and is a smooth hypersurface.
\item{(c)}
The ordinary locus consists in a smooth analytic curve in the deformation
space.
\item{(d)}
The inertially ample locus is equal to the complement of the union of
three hypersurfaces, any two of which meet transversely.

They also show that generic $S_3$-extensions actually exist:

\prop{\degenexist} ([B-M], prop.\ 9) Let $a$ be an integer such that
$27+4a^3$ is positive, prime, and less than $10^4$.  Then the splitting
field of the polynomial $x^3+ax+1$ is a generic $S_3$-extension for the prime
$27+4a^3$.  (I have verified this for $27+4a^3 < 10^{15}/2$, and believe that
there are no counterexamples.)

We will prove that the statements regarding the inertially reducible locus
and the ordinary locus still hold in the degenerate case, and we will say
something about the other loci as well.  In addition, we will give two
examples (of many) of degenerate $S_3$-extensions, one of them the splitting
field of the polynomial $x^3+7x-12$.

I would like to thank Barry Mazur and Nigel Boston for discussions of this
problem, and Fernando Gouv\^ea for encouraging me to pursue it after a long
hiatus.  Calculations in this paper were done using {\tt gp} and the tables of
number fields prepared at the Universit\'e de Bordeaux.

\section{Basics}

We start with some fundamental
definitions borrowed from [B-M], with very slight modifications.

\define (Cf. [B-M], Definition 2.)  Let $L/\Q$ be a totally complex
$S_3$-extension in which $p$ splits as
$(\p_1\p_2\p_3)^2$, and let $S$ be the set of finite ramified primes of $L$.
We say that $L$ is {\sl admissible} for $p$, or $L$ is {\sl neat}, if:

\item{1.} Any global unit of $L$ which is locally a $p$\/th power at all
elements of $S$ is globally a $p$\/th power.
\item{2.} The class number of $L$ is prime to $p$.
\item{3.} The completion of $L$ at any element of $S$ does not contain $p$\/th
roots of $1$.  (In particular, it follows that the cardinality of the residue
field is not congruent to $1$ mod $p$.)
Let $L$ be a neat $S_3$-extension of $\Q$, and let $\p_1, \p_2, \p_3$
be the primes of $L$ lying above $p$.  Let $e_1, e_2$ be a basis for global
units mod $p$\/th powers.  Since we are assuming that $L$ is neat, we may
suppose that $e_1$ is not a $p$\/th power in $L_{\p_1}$, and we may also
arrange things so that $e_2$ is not a $p$\/th power in $L_{\p_2}$ or
$L_{\p_3}$.
If $e_2$ is not a $p$\/th power in $L_{\p_1}$ either, then $L$ is generic.

\define The {\sl degeneracy index} of $L$ at $p$ will be the largest integer
$i$ such that $e_2$ is a $p^i$\/th power in $L_{\p_1}$.
(Of course this is a finite number, for the only elements of $L_{\p_1}$ which
are $p^i$\/th powers for all $i$ are $p-1$-st roots of $1$.)

The authors of [B-M] pay particular attention to the Galois closures of
cubic fields of the form $\Q(x)$, where $x^3+ax+1 = 0$, for $a$ an integer
such that $27+4a^3$ is positive and prime.  They show that the first seven
such fields are generic $S_3$-extensions of $\Q$, using a simple numerical
criterion.  As noted above, I have extended this verification to all
$a < 500000$, and find it hard to believe that there are any counterexamples.
However, if one does not restrict to these particular cubic fields, it becomes
easy to find degenerate $S_3$-extensions.  (The tables of number fields
available by anonymous FTP from {\tt megrez.math.u-bordeaux.fr} greatly 
facilitate such a search.)

The rest of the paper will be devoted to modifying the proofs and results of
Boston and Mazur so that they apply in the degenerate case.  That is,
we will determine the natural subspaces of the universal deformation space,
just as they do in their generic situation.  Regrettably, I do not have
anything to say about the cases which are not neat, whether because the
class number is a multiple of $p$ or because there is a unit which is a
$p$\/th power locally at all primes over $p$ but not globally.

\section{Definitions and Notations}
We now recall some more definitions from [B-M].

\define Let $L$ be an $S_3$-extension of $\Q$, and let $p$ be a rational
prime greater than $3$ which decomposes in $L$ as $\p_1\p_2\p_3$.  (We assume
that such a prime exists.)  Let $S$ be the
set of ramified primes of $L$.  Let $P$ be the Galois group over $L$ 
of the maximal
pro-$p$ extension of $L$ unramified away from $p$, or outside $S$ (in the
situations we will be considering, these are the same), $G$ its Galois
group over $\Q$, $L_p$ the completion of $L$ at $\p_1$, $P_p$ the Galois group
over $L_p$ of its maximal pro-$p$ extension, and $G_p$ the Galois group
of the maximal pro-$p$ extension of $L_p$ over $\Q_p$.  We also fix an
embedding of $\bar \Q$ into $\bar \Q_p$, and thus of $\Gal L_p$ into $\Gal L$,
such that the inertia subgroup $P^0_p$ maps to the inertia subgroup for $\p_1$.

\prop{\galoisgps}
$P$ is a free pro-$p$ group on $4$ generators, and $P_p$ is a free
pro-$p$ group on $3$ generators.

\proof [B-M, props.\ 3, 4]. \qed

We take $\sigma$ (resp.\ $\tau$) to be an element of order $2$ (resp.\ $3$) in
$S_3$.  Following one of the notations in [B-M], we will let $P$ be
generated by $u, \tau(u), \tau^2(u), v$, where $u$ conjugated by $\tau$ is,
obviously, $\tau(u)$, 
$\tau(v) = v$, $\sigma(u) = u$, and $\sigma(v) = v^{-1}$.  On the other hand,
$P_p$ will be generated by $\xi, \eta, \phi$, with $\xi$ and $\eta$ generating
the inertia and the nontrivial element of $\Z/2\Z$ acting as $+1$ on $\xi,
\phi$ and $-1$ on $\eta$.

\define Let $E$ be the group of global units of $L$.  It is the direct
product of a free abelian group of rank $2$ with a cyclic group of order $2$.
For any place $v$ of $L$, let $E_v$ be the group of units in the ring of
integers of $L_v$.

\define For any topological group $T$, let $\bar T$ be its $p$-Frattini
quotient.  More generally, let $\subbar{i} T$ be the maximal quotient of
$T$ which is an abelian pro-$p$ group with exponent dividing $p^i$ (that is,
$\subbar{i} T = T/(T,T)T^{p^i}$).

\define Let $K$ be a cubic extension of $\Q$.  Let $L$ be its Galois closure,
and let $S$ be the set of finite ramified primes of $L$.
Global class field theory gives us a map from the id\`ele class group of $L$
to the abelianization of its absolute Galois group.  This induces a map
$\oplus_{v \in S} \bar E_v \mapsto \bar P$ which is trivial on the image in
$\oplus_{v \in S} \bar E_v$ of the global units.  Under the conditions that
the class number of $L$ be prime to $p$ and that no completion of $L$ at
a prime in $S$ contain the $p$\/th roots of $1$, this map is surjective.

In this situation we say, as above, that $L$ is {\sl neat} for $p$, or for
$S$, if the map $\bar E \mapsto \oplus_{v \in S} \bar E_v$ is injective.
In this case, we consider a map $\bar E \mapsto \bar E_1$.  If it too is
injective, we are in the {\sl generic} situation treated by Boston and Mazur.
Otherwise, the extension is termed {\sl degenerate}, as remarked above, and
the degeneracy index is the largest $i$ for which the map $\subbar{i} E \mapsto
\subbar{i} E_1$ has cyclic image.

At this point we give our promised example of a degenerate $S_3$-extension.

\example{\degen} Let $K$ be the field $\Q(x)$,
where $x^3+7x-12 = 0$, and $L$ its Galois closure.
We claim that $L$ is a degenerate $S_3$-extension of
$\Q$, of degeneracy index $1$, with $p=5$.  Here the set of ramified primes
is $\{5,263\}$, both of which split as $(\p_1\p_2\p_3)^2$, so it is easy to check
that the completions there do not contain fifth roots of $1$.  The class number
of $L$ is $2$.

Using {\tt gp}, it is easy to check that the units of the cubic subfields of
$L$ generate the full unit group of $K$.  A fundamental unit of $K$ is
$14x-19$.  $K$ has a unique embedding into $\Q_5$, in which the image of $x$ is
congruent to $62$ mod $125$, so that the image of the fundamental unit is
congruent to $-1$ mod $25$, but not mod $125$, and is therefore a fifth power
but not a $25$th power.  I assert that the image of $x$ under
the embedding of $K$ into $\Q_5(\sqrt{10})$ is not a fifth power.  This will
essentially complete the verification that $L$ is neat.

In fact, it is easy to show that a unit of $\Q_5(\sqrt{10})$ which is
congruent to $1$ modulo $\m$, the maximal ideal, is a fifth power iff it
is congruent to $1$ modulo $\m^3$.  However, a root of $x^3+7x-12$ which
does not belong to $\Q_5$ is congruent modulo $\m^2$ to $4 \pm \sqrt{10}$, and
thus the unit is congruent modulo $\m^2$ to $2 \mp \sqrt{10}$
and cannot be a fifth power (multiply by $3^5$).

So, if we take a unit $u = \pm u_1^a u_2^b$ of $L$, where $u_1, u_2$ are
fundamental units in different cubic subfields of $L$, then in one completion
$u$ is a fifth power iff $5 | a$, and in another iff $5 | b$.  Thus $u$ is
a fifth power locally iff it is a fifth power globally, which completes the
proof that $L$ is neat.

\example{\highindex}
Let $p$ be a prime, $n$ a positive integer, and $r$ and $s$ integers
with $p^n | (r+s)$.  Suppose that the polynomial $x^3+rx^2+sx-1$ is irreducible
and does not have all real roots, and let $K$ be its root field and $L$ its
splitting field as above.  Suppose further that $p$ decomposes in $K$ as
$\p_1^2\p_2$.  (This cannot be arranged for all $p$, since it requires the
polynomial $x^4-8x^3+18x^2-27$ to have a root mod $p$, but it can for some, for
example with $p = 5, r = -29, s = 4, n = 2$.)  Clearly $K$ has a unit
congruent to $1$ modulo $p^n$ in the embedding of $K$ into $\Q_p$.
This must be the $p^{n-1}$-th power of a local unit,
and it seems that it is always possible to
arrange for $L$ to be neat by varying $r, s$ within their congruence
classes mod $p^n$.  Constructions like this lead me to
believe that for all primes $p > 3$, there are $S_3$-extensions with all
degeneracy indices at $p$, but I cannot prove it.  Among other problems, it can
never be possible to estimate the class number of $K$ as being less than $p$,
as did Mazur [M, section 1.13], because for a fixed $p$ the discriminant must
grow as $n$ increases.

\section{Representation Theory}

Let $G$ be a finite group and $F$ a field of characteristic prime to $\card G$.
Then, of course, the group algebra $F[G]$ is semisimple.  It is isomorphic
to a direct sum of matrix algebras over $F$ iff all irreducible representations
of $G$ (say there are $c$) can be defined over $F$.

Suppose we are in this case, and let $R$ be a local Artinian ring with residue
field $F$.  Since $R[G] \isisom F[G] \otimes_F R$, it is clear that $R[G]$ is
likewise a direct sum of $c$ matrix algebras.  Now, representations of $G$
with coefficients in $R$ correspond naturally to $R[G]$-modules free over $R$.
These, then, correspond to $R^c$-modules free over $R$,
that is, to $c$-tuples of free $R$-modules.  In turn, these correspond
canonically to $c$-tuples of $F$-modules, whence to $F[G]$-modules or to
representations.  In summary, a representation of $G$ to $M_n(R)$ is uniquely
determined up to conjugacy by its reduction to $M_n(F)$.  The usual theorems
on reducibility of representations then follow for representations to $R$.
For example, if we have a representation $\rho$ to $R$ which is an extension
of representations, it must in fact be their sum, for $\rho$ and the sum
have the same reduction to $M_n(F)$.

We will apply these ideas with $F = \F_p$ and $G = S_3$.  In essence, they
allow us to immediately take over all results about module decompositions
given in [B-M] without change here.  Since $P$, the Galois group of the
maximal pro-$p$ extension of $L$ unramified away from $p$ over $L$, is a
free pro-$p$ group, the groups $\subbar{i} P$ are all free
modules over $\Z/p^i\Z$ of the same rank.

\prop{\actions} $A$ is a semidirect product of $P$ by $S_3$ such that for
any $j$, $A$ acts
on $\subbar{j} P$ by $V + \chi$, where $\chi$ is the
nontrivial $1$-dimensional representation of $S_3$ with
coefficients in $\F_p$ and $V$ is the natural $3$-dimensional representation
of $S_3$.  For any $j$, $A_p$ is a semidirect product of $P_p$ by $\Z/2\Z$,
with $\Z/2\Z$ acting on $\subbar{j} P_p$ by $1 + 1 + \chi$.  In the global
case, the inertia subgroup maps to the space spanned by a basis vector of $V$
and $\chi$ (not a submodule, since different choices of prime above
$p$ give different inertia subgroups); in the local case, to $1 + \chi$.

\proof [B-M, props.\ 7 and 8], together with the above to remove the
restriction $j = 1$ made there. \qed

Boston and Mazur study the exact sequence of $p$-Frattini quotients
$$0 \mapsto \bar E \mapsto \bar E_1 \oplus \bar E_2 \oplus \bar E_3 \mapsto
\bar P \mapsto 0.$$  Likewise we will study the exact sequence of
$p^i$-quotients.  That is, we define a map $\Pi_j^i$ to be that given
by class field theory from $\subbar{j} E_k$ to $\subbar{j} P$.  Its image will
be denoted $\subbar{j} P_k$, and by class field theory $\subbar{j} P_1$ is
the image of the inertia subgroup $\subbar{j} P^0_p$ in $\subbar{j} P$.

\prop{\inter} Let $L$ be an $S_3$-extension of $\Q$, degenerate for $p$ with
degeneracy index $i$.  Then the intersection of any two, or all three, of
the $\subbar{j} P_k$ is
isomorphic to $\Z/p^l\Z$, where $l = \min(i,j)$.  (In the case $j = 1$, this
reduces to the results in the first part of [B-M, section 2.3].)

\proof We consider the cases $i < j, \, j \le i$ separately.  In the case
$j \le i$, the image of $P_1$ is isomorphic to $\Z/p^j\Z$, and it is
stable under the action of the involution of the Galois group which fixes
$\p_1$.  The rest of the proof in [B-M, section 2.3] can now be taken over
word for word. \qed

We now consider the case $i < j$.  Everything is compatible with the inclusion
maps $\subbar{j} E_k \mapsto \subbar{j+1} E_k$ and $\subbar{j} P \mapsto
\subbar{j+1} P$, so the image must contain a $\Z/p^i\Z$-subgroup, and no more
elements of order dividing $p^i$.  If there is an element of higher order
in the intersection $\subbar{j} P_k \cap \subbar{j} P_{k'}$, say $y$, coming
from $y_k$ and $y_{k'}$, then the element $y_k \oplus y_{k'} \oplus 0$ would
be in the image of $\subbar{j} E$, by exactness.  This would immediately imply
that the degeneracy index of $L$ is greater than $i$.\qed

\section{Linking Local and Global Presentations}

We have already described (in \actions) the presentations of the local and
global Galois groups $G_p, G$.  Now we must show how they behave under the
map $G_p \mapsto G$ (in particular, what happens when we restrict this to
a map $\Pi_p \mapsto \Pi$).  This is where the difference between the generic
and degenerate situations becomes important.

\prop{\module} (Cf.\ [B-M, lemma 2.4.4].)  Suppose that $i$, the degeneracy
index of $L$, is at least $j$, and let $\xi, \eta$ be generators of the inertia
subgroup of $\subbar{j} P_p$ such that the nontrivial element of $\Gal(L_p/\Q)$
acts as $+1$ on $\xi$ and $-1$ on $\eta$.  Let $r, s$ be the images of $\xi, \eta$ in
$\subbar{j} \Pi$, and let $R, S$ be the $S_3$-stable subspaces that they
generate.  Then $R \isisom 1 + \epsilon$ and $S \isisom \chi$.

\proof Recall that $\subbar{j} \Pi \isisom 1 + \chi + \epsilon$.  Because
$L$ has no unramified extensions of degree $p$, the $S_3$-stable subspace
of $\subbar{j} \Pi$ generated by the image of a local inertia group---that
is, $R + S$---must be
the whole thing.  Also, $R, S$ must be quotients of the inductions of
$1$ and $\chi$ from $A_p = \Z/2\Z$ to $A = S_3$, respectively.  On the other
hand, if $i < j$, neither $R$ nor $S$ can be one-dimensional over
$\Z/p^j\Z$, by \inter.  Thus, the
statement on $R$ follows if we prove the statement about $S$.

Corresponding to $\xi, \eta$, we let $a, b$ be generators of $\subbar{j} E_1$
such that $a^\sigma = a, b^\sigma = 1/b$.  Let $\upsilon$ be the global
unit of $L$ which is a $p^j$\/th power in $L_{\p_1}$.  Then the
image of $\upsilon$ in $\subbar{j} E_{\p_2}$ must be a multiple of $b$.
Indeed, on the one hand, the product of the three global conjugates of
$\upsilon$ can be taken to be $1$, and on the other hand, the two conjugates
that are not in $\Q_p$ are local conjugates in $L_p$, so when reduced to
$\subbar{j} E_{L_p}$ they have the same coefficient of $a$, which must
therefore be $0$.  On the other hand, the coefficient of $b$ must be a unit,
for otherwise $\upsilon$ would be a $p$\/th power everywhere locally, a
possibility excluded by our hypotheses.

In particular, the image of $\subbar{j} E_L$ in $\oplus_k\, \subbar{j} E_k$
is spanned by $(b,-b,0)$ and $(0,b,-b)$.  It follows that the intersection
of the images of the $\subbar{j} E_k$ in $\subbar{j} P$ is the image of
$(b,0,0)$, which is obviously the $\chi$ subspace as claimed.\qed

A curious consequence of this proposition is as follows:

\cor{\classno} Let $L/\Q$ be an $S_3$-extension of degeneracy index at least
$j$ for $p$, and let $F$ be the subextension of $\Q(\zeta_{p^{j+1}})$ which
is of degree $p^j$ over $\Q$.  Then the class group of the compositum $L \vee
F$ contains a subgroup isomorphic to $(\Z/p^j\Z)^2$.

\proof It is sufficient, of course, to construct an unramified extension with
this Galois group.  The point is simply that the inertia groups for
$\p_i$ in the extension cut out by $\subbar{j} P$ are isomorphic to
$(\Z/p^j\Z)^2$ and have an intersection isomorphic to $(\Z/p^j\Z)$, which cuts
out the extension $M$, say.  On the other hand, $L \vee F
\subset M$ and is totally ramified over $L$ at each $\p_i$.  Thus,
$M/(L \vee F)$ is unramified at these primes, and (by definition
of $P$) at all others as well. \qed

We must now specify the relation between local and global presentations
more precisely.

\prop{\link} (Cf.\ [B-M, prop.\ 10].)
Let $L$ be an admissible $S_3$-extension of $\Q$ of degeneracy
index $i$ for the prime $p$.  Then we may take the local and global
systems of generators such that the image of $\xi$ is $u$ and,
in the induced map on quotients $\subbar{j} \Pi_p \mapsto \subbar{j} \Pi$, 
the image of $\eta$ is $v$, if $i \ge j$.
\proof By [B-M, prop.\ 7 and addendum] we may take the image of $\xi$ to
be $u$.  The statement about $\eta$ follows from the last lemma, similarly to
the proof of [B-M, prop.\ 10]. \qed

We have now accumulated all necessary information about the Galois groups and
can proceed to studying the universal deformation.

\section{The Universal Deformation}

Let $L$ be an admissible $S_3$-extension; for the moment, the index of
degeneracy does not matter.  There is a Galois representation
$\bar \rho: G \mapsto GL_2(\F_p)$, unique up to conjugacy, which factors
through $\Gal L/\Q$ and maps it injectively into $GL_2(\F_p)$.
For concreteness, we fix elements $\sigma$, $\tau$
in $S_3$ of order $2$ and $3$ respectively and map them to 
$$\pmatrix{1&0\cr 0&-1\cr}, \pmatrix{-1/2& 1/2 \cr -3/2 & -1/2\cr}.$$

We will be studying deformations of $\bar \rho$ to complete local noetherian
rings with residue field $\F_p$.  The universal deformation has been completely
described.

\prop{\univdef} The universal deformation ring is
the power series ring $\Z_p[[T_1,T_2,T_3]]$, and the universal deformation
may be given as follows:
$$\sigma \mapsto \pmatrix{1&0\cr 0&-1\cr},
\tau \mapsto \pmatrix{-1/2&1/2\cr -3/2&-1/2\cr},
u \mapsto \pmatrix{1+T_1&0\cr 0&1+T_1\cr},
v \mapsto \pmatrix{(1-3T_3^2)^{1/2}&T_3\cr -3T_3&(1-3T_3^2)^{1/2}\cr}.$$

\proof This is [B-M], prop.\ 11, and a detailed proof is given there. \qed

To understand the universal deformation more fully, we must understand the
image of $\eta$.  Since $\eta$ conjugated by $\sigma$ is $\eta^{-1}$, the
image of $\eta$ must have determinant $1$ and equal diagonal entries, so it
is, say, $$\pmatrix{(1+fg)^{1/2}&f\cr g & (1+fg)^{1/2}\cr}.$$

\prop{\imeta} Modulo $\m$, the power series $f$ is congruent to $T_3$, and
$g$ to $-3T_3$.

\proof As [B-M], prop.\ 12, except that here the image of $\subbar{1} 
\eta$ under the natural map $\bar \Pi_p \mapsto \bar \Pi$ is $\bar v$.\qed

We can now determine some of the natural subspaces.  We will be considering
representations of the Galois group into $GL_2(\Z_p)$ which are deformations
of the representation into $GL_2(\F_p)$.  Thus they come from the universal
deformation, and are described by a continuous homomorphism
$\Z_p[[T_1,T_2,T_3]] \mapsto \Z_p$.  Such a homomorphism $\alpha$ is described
by giving $\alpha(T_1), \alpha(T_2), \alpha(T_3)$; the space of such is
therefore naturally identified with $p\Z_p \times p\Z_p \times p\Z_p$, which
is a $3$-dimensional $p$-adic manifold.
The only visible difference between our situation and the generic one is that
here $f$ and $g$ are not transversal.  Presumably the order of contact of
their zero loci is equal to the degeneracy locus of the extension, but I do
not see how to prove this.

\prop{\inred}
(cf. [B-M], prop.\ 13.)  The inertially reducible locus is the union
of the hypersurfaces $f = 0$ and $g = 0$.  The ordinary locus is the smooth
curve defined by $T_1 = g = 0$.

\proof Identical to the proofs given in [B-M]. \qed

It is still true that a representation is inertially dihedral iff $T_1 = T_2$
or $f = g = 0$.  If, as is presumably the case, the loci $f = 0$ and $g = 0$
are distinct, the argument in [B-M] goes through to show that $f = g = 0$
implies $T_1 = T_2$.  (This sounds like something that should be easy to
prove, but I have not managed to.)  It would then follow, just as in [B-M],
that the inertially ample locus is the complement of the union of the
inertially reducible and inertially dihedral loci.

\secnum = -1
\section{References}
\frenchspacing\leftskip=20.0pt \parindent=-20.0pt\parskip=10pt
[B-M] N. Boston, B. Mazur, {\sl Explicit universal deformations of
Galois representations}.  In {\sl Algebraic Number Theory}, Adv. Stud. Pure
Math. {\bf 17}, 1--21.

[M] B. Mazur, {\sl Deforming Galois representations}.  In {\sl Galois groups
over $\Q$}, MSRI Publications {\bf 16}, 385--437.

\bye